\documentclass[10pt]{amsart}
\usepackage{amsfonts}
\usepackage{graphicx}
\usepackage{multirow}
\usepackage{url}
\usepackage[colorlinks,linkcolor=blue,citecolor=red]{hyperref}

\theoremstyle{plain}
\newtheorem{Thm}{Theorem}

\newtheorem{Lem}{Lemma}
\theoremstyle{definition}

\renewcommand{\Re}{\operatorname{Re}}

\allowdisplaybreaks[3]

\title[Calculate the Riemann zeta function]{Numerical calculation of the Riemann zeta function at odd integer arguments: A direct formula method}

\author[Q.~Luo]{Qiang Luo}
\address{Department of Physics, Renmin University of China, Beijing, 100872, P. R. China}
\email{qiangluo@ruc.edu.cn}

\author[Z.~D.~Wang]{Zhidan Wang}
\address{School of Mathematical Science and Technology, Yangzhou University, Yangzhou, Jiangsu 225002, P. R. China}
\email{zhitanwang@gmail.com}

\begin{document}
\begin{abstract}
In this article, we introduce a recurrence formula which only involves two adjacent values of the Riemann zeta function at integer arguments. Based
on the formula, an algorithm to evaluate $\zeta$-values(i.e. the values of Riemann zeta function) at odd-integers from the two nearest $\zeta$-values at even-integers is posed and proved. The behavior of the error bound is $O(10^{-n})$ approximately where $n$ is the argument. Our method is especially powerful for the calculation of Riemann zeta function at large argument, while for smaller ones it can also reach spectacular accuracies such as more than ten decimal places.
\end{abstract}
\maketitle
\section*{Introduction}
\setlength\parindent{0em}
Zeta functions of various kinds, such as Hurwitz zeta function, Epstein zeta function and Dirichlet $L$-function, are all-pervasive objects in modern mathematics, especially in analytical number theory, and among which the prototype zeta function is the famous Riemann zeta function. It is classically defined as the sum of the infinite series\cite{ref1,ref2,refEdwards}
\begin{equation}\label{eq:1Riemann}
\zeta(s)=\sum_{n=1}^{\infty}\frac{1}{n^s}
\end{equation}
with the complex variable \begin{math}s=\sigma+it\end{math}. Specially, the series converges if \begin{math}\sigma=\Re s>1\end{math}. We can extend $\zeta(s)$ from $s$ with $\Re s > 1$ to $s$ with \begin{math}\Re(s)>0, s\not=1\end{math} by the following formula
\begin{equation}\label{eq:2end}
\eta(s)=\sum_{n=1}^{\infty}\frac{(-1)^{n+1}}{n^s}=(1-2^{1-s})\zeta(s)
\end{equation}
where $\eta(s)$ is the Dirichlet eta function or alternating eta function.

Historically, people prefer to study the closed form of the Riemann zeta function at positive integer arguments in that those special values seem to dictate the properties of the objects they associated. In condensed matter physics for instance, the famous Sommerfeld expansion, which is usful for the calculation of particle number and internal energy of electrons, involves Riemann zeta function at even integers\cite{refSMFexpan}, while the spin-spin correlation functions of isotropic spin-$1/2$ Heisenberg model are expressed by $\ln2$ and Riemann zeta functions with odd integer arguments\cite{refHeisenberg}. It was, without doubt, a profound discovery of Euler in 1736 to work out the prolonged Basel problem\cite{ref3}
\begin{equation}\label{eq:1Eulersum}
\zeta(2)=\frac{{\pi}^2}{6}
\end{equation}
superbly. It is well-known that for positive even integer arguments the Riemann zeta function can be expressed explicitly as\cite{ref6}
\begin{equation}\label{eq:1Beven}
\zeta(2n)=\frac{(-1)^{n+1}(2\pi)^{2n}}{2(2n)!}B_{2n}
\end{equation}
in terms of the Bernoulli numbers \begin{math}B_{n}\end{math}. On the contrary, however, the explicit formula for Riemann zeta function at odd values is difficult if not fundamentally impossible to obtain. Euler himself once conjectured that $\zeta(2n+1)=c(n)\pi^{2n+1}$ and $c$ involves the irrational constant $\eta(1) = \ln 2$\cite{zeta3error}. This suggests that Riemann zeta function at odd integers produces a recurrence relation that is self-recursive. Even up to now, for positive odd integer arguments the Riemann zeta function can only be expressed by series and integral(see \eqref{zetafuncInt} and \eqref{eqZetaSeries} for detail). One possible integral expression is\cite{ref7}
\begin{equation}\label{eq:1Bodd}
\zeta(2n+1)=\frac{(-1)^{n+1}(2\pi)^{2n+1}}{2(2n+1)!}\int_{0}^{1}B_{2n+1}(x)\cot({\pi}x)dx
\end{equation}
where \begin{math}B_{2n+1}(x)\end{math} are Bernoulli polynomials. A relevant aspect is that, for Riemann zeta function, the celebrated Goldbach-Euler theorem\cite{refamm} assumes the elegant form
\begin{equation}\label{eq:1frac}
\sum_{n=2}^{\infty}\textbf{frac}(\zeta(n))=1,
\end{equation}
where $\textbf{frac}(x)=x-[x]$ denotes the fractional part of the real number $x$. It turns out that
\begin{equation}\label{eq:1fracoddeven}
\sum_{n=1}^{\infty}\textbf{frac}(\zeta(2n))=\frac{3}{4}, \sum_{n=1}^{\infty}\textbf{frac}(\zeta(2n+1))=\frac{1}{4}.
\end{equation}

Indeed, the formulas \eqref{eq:1Beven} and \eqref{eq:1Bodd}, along with \eqref{eq:1fracoddeven} do reveal somewhat similarity for the values of Riemann zeta function at even and odd arguments. Meanwhile, the calculation of Riemann zeta function and related series is a hot topic in computational mathematics. The traditional methods are Euler-Maclaurin formula and Riemann-Siegel formula, and algorithms are still being developed in earnest ever since\cite{integral,series1,series2,refLima}. Typically, a particular numerical method is limited to a special domain. Therefore, when concentrating on Riemann zeta function at odd integers, a special method should be constructed in view of the connection of Riemann zeta function values between odd and even integers.

In this paper we mainly obtain a recurrence formula \eqref{eq:4recur} relating to the Riemann zeta function and based on which we construct an algorithm for the calculation of the Riemann zeta function at odd integers. In addition, numerical calculation implies that the algorithm can reach considerable accuracies with small odd integer arguments, not to speak of larger ones. Quantificationally, the behavior of the error bound is $O(10^{-n})$ where $n$ is the argument.
\section*{Notations and Preliminaries}
We begin by recalling the definition of the Bernoulli polynomials \begin{math}B_{n}(x)\end{math} and their basic properties in a nutshell to render the paper essentially self-contained. The generating function of the Bernoulli polynomials \begin{math}B_{n}(x)\end{math} is \cite{ref1,ref2,refEdwards}
\begin{equation}\label{eq:2Ber_poly}
\frac{te^{tx}}{e^t-1}=\sum_{n=0}^{\infty}\frac{B_{n}(x)}{n!}t^n.
\end{equation}
Taking a derivative with respect to $x$ on both sides of \eqref{eq:2Ber_poly}, we find that
\begin{equation}\label{eq:2Bdiff}
B_{n}^{\prime}(x)=nB_{n-1}(x).
\end{equation}
Bernoulli polynomials can also be expressed explicitly from Bernoulli numbers
\begin{equation}\label{eq:2Num_Poly}
B_{n}(x)=\sum_{k=0}^{n}\left(\begin{array}{c}n\\k\end{array}\right)B_{k}x^{n-k}.
\end{equation}

For convenience, we introduce two kinds of reduced Bernoulli numbers(RBNs), one relates to the even-labeled Bernoulli numbers (denoted by $+$)
\begin{equation}\label{eq-2.5}
B_{n}^{+}=(-1)^{n+1}B_{2n},
\end{equation}
and another relates to the odd-labeled Bernoulli polynomials (denoted by $-$)
\begin{equation}\label{eq-2.6}
B_{n}^{-}=(-1)^{n+1}\int_{0}^{1}B_{2n+1}(x)\cot({\pi}x)dx.
\end{equation}
In this section we will demonstrate the asymptotic representation of the two kinds of RBNs in a uniform framework and establish their integral representation subsequently.
\subsection*{Asymptotic representations of RBNs}
The asymptotic expressions of Bernoulli polynomials at even and odd subscript are respectively\cite{ref11}
\begin{subequations}\label{eq3Bsim}
\begin{equation}\label{eq:3Bsima}
(-1)^{n+1}\frac{(2\pi)^{2n}}{2(2n)!}B_{2n}(x){\sim}\cos(2{\pi}x),
\end{equation}
\begin{equation}\label{eq:3Bsimb}
(-1)^{n+1}\frac{(2\pi)^{2n+1}}{2(2n+1)!}B_{2n+1}(x){\sim}\sin(2{\pi}x).
\end{equation}
\end{subequations}
Moreover, the Bernoulli polynomials can also be expressed in a stronger form based on Fourier sine and cosine series expansion\cite{ref12}
\begin{subequations}\label{eq3Bfourier}
\begin{equation}\label{eq:3Bfouriera}
B_{2n}(x)=(-1)^{n+1}\frac{2(2n)!}{(2\pi)^{2n}}\sum_{k=1}^{\infty}\frac{\cos(2{\pi}kx)}{k^{2n}},
\end{equation}
\begin{equation}\label{eq:3Bfourierb}
B_{2n+1}(x)=(-1)^{n+1}\frac{2(2n+1)!}{(2\pi)^{2n+1}}\sum_{k=1}^{\infty}\frac{\sin(2{\pi}kx)}{k^{2n+1}}.
\end{equation}
\end{subequations}
Obviously \eqref{eq3Bsim} is just the corollary of \eqref{eq3Bfourier}. Joining together we therefore obtain the asymptotic behavior of RBNs
\begin{subequations}
\begin{equation}\label{eq:3asy_Beven}
B_{n}^{+}=(-1)^{n+1}B_{2n}(0){\sim}\frac{2(2n)!}{(2\pi)^{2n}}
\end{equation}
\begin{equation}\label{eq:3asy_Bodd}
B_{n}^{-}=(-1)^{n+1}\int_{0}^{1}B_{2n+1}(x)\cot({\pi}x)dx{\sim}\frac{2(2n+1)!}{(2\pi)^{2n+1}}.
\end{equation}
\end{subequations}
where we use the fact that $\int_{0}^{1}\sin(2{\pi}x)\cot({\pi}x)dx=1$.
\subsection*{Integral representations of RBNs}
Let us consider two auxiliary integrals\cite{ref6}
\begin{subequations}
\begin{equation}\label{eq:3.7a}
I_{c}(n,m)=\int_{0}^{1}B_{2n}(t)\cos(m{\pi}t)dt
\end{equation}
\begin{equation}\label{eq:3.7b}
I_{s}(n,m)=\int_{0}^{1}B_{2n+1}(t)\sin(m{\pi}t)dt
\end{equation}
\end{subequations}
with $m$ and $n$ are integers and \begin{math}n\geq1\end{math}. Specially, when \begin{math}n=1\end{math}, direct computation shows that
\begin{subequations}\label{eqIntsinecosine}
\begin{align}\label{eq:3Ic1}
I_{c}(1,m)&=\int_{0}^{1}\Big(t^2-t+\frac16\Big)\cos(m{\pi}t)dt \nonumber\\
&=\left\{
  \begin{array}{cl}
  0                  &, m=1,3,5,\cdots\\
  \frac{2!}{(m\pi)^2}&, m=2,4,6,\cdots
  \end{array}
  ,
\right.
\end{align}
\begin{align}\label{eq:3Is1}
I_{s}(1,m)&=\int_{0}^{1}\Big(t^3-\frac32t^2+\frac12t\Big)\sin(m{\pi}t)dt \nonumber\\
&=\left\{
  \begin{array}{cl}
  0                   &, m=1,3,5,\cdots\\
  \frac{3!}{(m\pi)^3} &, m=2,4,6,\cdots
  \end{array}
  .
\right.
\end{align}
\end{subequations}
By virtue of \eqref{eq:2Bdiff} and integrating by parts twice, readily yields
\begin{subequations}\label{eqIcsrec}
\begin{equation}\label{eq:3Icrec}
I_{c}(n,m)=-\frac{(2n)(2n-1)}{(m{\pi})^2}I_{c}(n-1,m)
\end{equation}
\begin{equation}\label{eq:3Isrec}
I_{s}(n,m)=-\frac{(2n+1)(2n)}{(m{\pi})^2}I_{s}(n-1,m).
\end{equation}
\end{subequations}
Combining \eqref{eqIntsinecosine} and \eqref{eqIcsrec} we find that
\begin{subequations}
\begin{equation}\label{eq:3.10a}
I_{c}(n,m)=\frac{(-1)^{n+1}(2n)!}{(m{\pi})^{2n}}
\end{equation}
\begin{equation}\label{eq:3.10b}
I_{s}(n,m)=\frac{(-1)^{n+1}(2n+1)!}{(m{\pi})^{2n+1}}
\end{equation}
\end{subequations}
hold if $m$ is even. Immediately, the integral representations of the two RBNs \begin{math}B_{n}^{+}\end{math} and \begin{math}B_{n}^{-}\end{math} are
\begin{subequations}
\begin{equation}\label{eq:3.11a}
B_{n}^{+}{\sim}2(-1)^{n+1}I_{c}(n,2)
\end{equation}
\begin{equation}\label{eq:3.11b}
B_{n}^{-}{\sim}2(-1)^{n+1}I_{s}(n,2).
\end{equation}
\end{subequations}
\section*{Algorithm to calculate the Riemann zeta function}
Mathematically, Riemann zeta function is said to be monotonically decreasing since its values are only falling and never rising with increasing values of $s$ with $s\geq2$. Besides, $\zeta(2)=\frac{\pi^2}{6}$, $\zeta(+\infty)>1$, thus $0<\zeta(s)-1<1$. Analogously, one can show that $0<\frac{1}{\eta(s)}-1<1$. For brevity we denote the reciprocal function as below
\begin{equation}\label{eq:4def}
\rho(s)\equiv\textbf{frac}\Big(\frac{1}{\eta(s)}\Big)=\frac{1}{\eta(s)}-1.
\end{equation}
Now that what we concerned most is the values of the Riemann zeta function at integers for the moment, the asymptotic behavior of the ratio of the reciprocal function \eqref{eq:4def} at odd integers and even integers interests us. Motivated by \eqref{eq:1frac} and \eqref{eq:1fracoddeven}, we manage to demonstrate a formula, the so-called \textit{recurrence formula}(not in a strict sense, though), on condition that the argument is a positive integer. Motivated by the recurrence, we manage to construct an algorithm to compute the Riemann zeta function.
\subsection*{Demonstration of the recurrence formula}
\begin{Thm}
\textit{If $n$ is a positive integer such that \begin{math}n \geq 1\end{math}, the recurrence relation holds
\begin{equation}\label{eq:4recur}
\lim_{n\rightarrow\infty}\frac{\rho(2n+1)}{\rho(2n)}=\frac12.
\end{equation}}
\end{Thm}
\begin{proof}
By using of \eqref{eq:1Beven} and \eqref{eq:1Bodd} and the definition of reciprocal function \eqref{eq:4def}, we have  %
\begin{equation}\label{eq:4.4}
\frac{\rho(2n+1)}{\rho(2n)}=\frac{(2n+1)!-(2^{2n}-1){\pi}^{2n+1}B_{n}^{-}}{(2n)!-(2^{2n-1}-1){\pi}^{2n}B_{n}^{+}}\frac{B_{n}^{+}}{{\pi}B_{n}^{-}}\frac{2^{2n-1}-1}{2^{2n}-1}.
\end{equation}
Since the limitation of the rightmost term is exactly equal to \begin{math}1/2\end{math} if $n$ is large enough, what we want to prove is
\begin{equation}\label{eq:4.5}
\frac{(2n+1)!-(2^{2n}-1){\pi}^{2n+1}B_{n}^{-}}{(2n)!-(2^{2n-1}-1){\pi}^{2n}B_{n}^{+}}\frac{B_{n}^{+}}{{\pi}B_{n}^{-}}{\sim}1
\end{equation}
or equivalently
\begin{equation}\label{eq:4.6}
\frac{(2n+1)!}{{\pi}^{2n+1}}\frac{1}{B_{n}^{-}}-\frac{(2n)!}{{\pi}^{2n}}\frac{1}{B_{n}^{+}}{\sim}2^{2n-1}.
\end{equation}
From the asymptotic formulae \eqref{eq:3asy_Beven} and \eqref{eq:3asy_Bodd} of the two kinds of RBNs we find that, without any difficulty, we have finished the demonstration of the recurrence formula of the Riemann zeta function.
\end{proof}
As a matter of fact, we can extend the validity of \eqref{eq:4recur} from positive integers to positive real numbers straightforward. We therefore can obtain the asymptotic behavior of Riemann zeta function as
\begin{equation}\label{zetaValueFromRecurrence}
\frac{1}{\zeta(s)}{\sim}\frac{2^{s-1}-1}{2^{2s-3}}\Big(\frac{2}{\zeta(2)}-1\Big)+\frac{2^{s-1}-1}{2^{s-1}}.
\end{equation}
by using of \eqref{eq:4recur}. The application of \eqref{zetaValueFromRecurrence} can be diverse, here we just pick a example relating to prime number theorem. The positive integer $x$ is $s$-free if and only if in the prime factorization of $x$, no prime number occurs more than $s$-1. Indeed, if \begin{math}Q(x,s)\end{math} denotes the number of $s$-free integers (e.g. 2-free integers being square-free integers) between 1 and $x$, one can show that\cite{refsFree}
\begin{equation}\label{eq:5Q}
Q(x,s)=\frac{x}{\zeta(s)}+O(\sqrt[s]{x}),
\end{equation}
therefore we find that the asymptotic density of $s$-free integers \begin{math}Q(x,s)/x{\sim}\frac{1}{\zeta(s)}\end{math} is nothing but \eqref{zetaValueFromRecurrence}.

\begin{figure}[!htb]  
\centerline{\includegraphics[width=8.0cm]{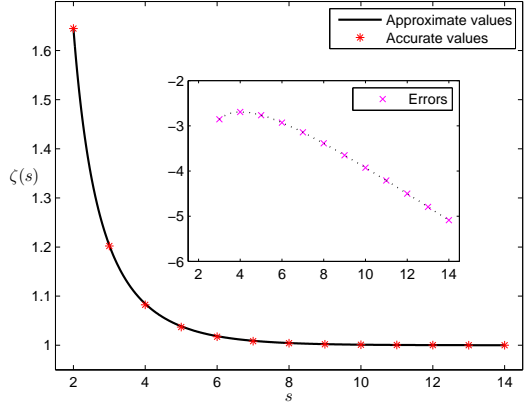}\hspace{4mm}}
\caption{\textit{Asymptotic behavior of Riemann zeta function}.~
The solid line represents the approximate values($\zeta^{\bf{ap}}(s)$) obtained by \eqref{zetaValueFromRecurrence}, while the stars "*" represent the accurate values($\zeta^{\bf{ac}}(s)$) when $s$ is an integer. The crosses "$\times$" in the inserted figure, indicate the base 10 logarithm of the absolute errors($\epsilon(s)=\lg\big(\vert{\zeta^{\bf{ap}}(s)-\zeta^{\bf{ac}}(s)}\vert\big)$) at integers.}\label{figError}
\end{figure}

Another intriguing issue is to what degree can \eqref{zetaValueFromRecurrence} reveals its ability to obtain the $\zeta$-values. Figure \ref{figError} is thus plotted as follow.

The fact that all the stars "*" lie on the solid curve indicates that \eqref{zetaValueFromRecurrence} may be a suitable candidate for the calculation of Riemann zeta function. The emergence of the abnormal slope between $s=3$ and $s=4$ in the inserted figure, however, implies that any $\zeta$-value obtained from its nearest neighbors should be much more accuracy. We therefore come up with a satisfactory proposal which is postponed until next subsection.
\subsection*{Basic ideas for the algorithm}
Abundant methods to evaluate the $\zeta(2n)$ have appeared in the mathematical literatures from now and then ever since Euler's seminal work. In contract, the explicit formula for odd-argument $\zeta$-values remains to be an open problem though some results shed light on it\cite{HeTX1,HeTX2}. By analogy to $\zeta(2n)$, several authors have established the series and integral representations of $\zeta(2n+1)$, which, to some degree, provides some perspectives on the difficulty of evaluating $\zeta(2n+1)$ as opposed to $\zeta(2n)$. From the viewpoint of numerical method, one natural way to construct the corresponding algorithm to evaluate the odd-argument Riemann zeta function is by viture of the even-argument $\zeta$-values near to them. In the current paper, only the two nearest $\zeta$-values are taken into consideration currently for simplicity. When $n$ is large enough, \eqref{eq:4recur} can be rewritten as
\begin{subequations}
\begin{equation}\label{eq:5up}
\rho^{l}(2n+1){\sim}\frac12\rho(2n)
\end{equation}
\begin{equation}\label{eq:5low}
\rho^{r}(2n+1){\sim}2\rho(2n+2)
\end{equation}
\end{subequations}
where $\rho^{l}(2n+1)$ and $\rho^{r}(2n+1)$ represent two different representations of the asymptotic behavior of $\rho(2n+1)$. Judging by appearance, One can use any of the formula above to calculate the Riemann zeta function at odd integers. When considering that those two formulae give the upper and lower bound of the zeta-values at odd integers(see Theorem 2), we come up with the idea that we can combine them together by a special method. It happens to us that there may exist a somewhat mysterious map from $\zeta(2n)$ and $\zeta(2n+2)$ to $\zeta(2n+1)$, which will ensure us to obtain the approximation values of $\zeta(2n+1)$ with higher precision. Let us give a proposition relating to Dirichlet eta function firstly before we move forward to give another theorem.
\begin{Lem}\label{Olss}
\textit{If $n$ is a positive integer such that \begin{math}n \geq 1\end{math}, the two inequalities hold
\begin{equation}\label{eq:5ineq1}
\frac{4}{\eta(2n+2)}-\frac{1}{\eta(2n)}>3
\end{equation}
\begin{equation}\label{eq:5ineq2}
\eta(2n)>\frac{2^{2n-1}-2}{2^{2n-1}-1}
\end{equation}}
\end{Lem}
Those two inequalities are quite new to the authors because we haven't seen them in any literature or monograph before. However, we are not intended to give the details here since the demonstration is rather elementary. The theorem below holds once we take advantage of lemma 1.
\begin{Thm}
\textit{If $n$ is a positive integer such that \begin{math}n \geq 1\end{math}, the inequality holds
\begin{equation}\label{eq:4ineqchain}
\zeta^{l}(2n+1)>\zeta^{r}(2n+1)>1
\end{equation}}
where $\zeta^{l}(2n+1)$ and $\zeta^{r}(2n+1)$ correspond to $\rho^{l}(2n+1)$ and $\rho^{r}(2n+1)$ respectively.
\end{Thm}

Since Riemann zeta function is a monotonic decreasing function, the exactly value $\zeta(2n+1)$ is just between $\zeta^{l}(2n+1)$ and $\zeta^{r}(2n+1)$ for any given positive integer $n$. For the benefit of accuracy we regard the geometric mean values of the \eqref{eq:5up} and \eqref{eq:5low} as the approximate values of the reciprocal function $\rho(2n+1)$, namely
\begin{equation}\label{eq:5mean}
\rho(2n+1){\approx}\sqrt{\rho(2n)\rho(2n+2)}
\end{equation}
which is the most valuable ingredient of our algorithm.

The basic steps for the calculation of $\zeta(2n+1)$ are presented as follow. Firstly, $\rho(2n)$ and $\rho(2n+2)$ should be calculated from \eqref{eq:1Beven}, \eqref{eq:2end} and \eqref{eq:4def}) in sequence. Secondly, the value of $\rho(2n+1)$ is ready to be obtained in light of \eqref{eq:5mean}. Lastly, the ultimate aim, i.e. $\zeta(2n+1)$ is just at hand from \eqref{eq:4def}) and \eqref{eq:2end}, reversely. Our algorithm doesn't bother circulation of any kind, it just looks like a formula, therefore we refer it as the \textit{direct formula method}.

In order to start our method, we need to know some $\zeta$-values at even integers. For example, $\zeta(2)$ and $\zeta(4)$ should be available to get $\zeta(3)$. We can obtain $\zeta(2n)$  through \eqref{eq:1Beven} systematically for small argument. However, it is almost impossible to obtain Bernoulli numbers by the ordinary recursive methods thus we hardly know the values of $\zeta(2n)$ if the argument is large enough. Many methods for computing Bernoulli numbers have been invented. David Harvey introduced an efficient multimodular algorithm\cite{bernmm} which ensures us to obtain the Bernoulli numbers $B_n$ at $n=10^8$. However, one can also use the intrinsic function \textbf{Zeta}[$s$] in Mathematica since it is also based on an efficient algorithm. Therefore, for convenience, our computation platform is mainly on Mathematica and we regard those values as benchmarks.
\section*{Calculation of Riemann zeta function at odd integers}
The calculation of Riemann zeta function plays an essential role in the study of number theory and associated subjects such as statistical physics and condensed matter physics. Various approaches to accomplish this task have been proposed\cite{integral,series1,series2,ref14}, especially for the evaluation of zeta function at integer arguments or in the critical strip (for the computation of Riemann's zeros). Most of the methods available consist of using integral forms of some particular functions or recursive series forms. Quite recently, Babolian \textit{et al} transform $\zeta(s)$ to some appropriate integral forms and introduce a method to compute the Riemann zeta function based on Gauss-Hermite and Gauss-Laguerre quadratures\cite{integral}. Numerical result show that 20 points are capable of producing an accuracy of seven-decimal place for small arguments. Besides, many rapidly converging series for $\zeta(2n+1)$ have been introduced by Srivastava in a review article\cite{series2} and by other authors\cite{series1,refLima}. In this section we firstly give some numerical examples according to our method to illustrate its accuracy, then we compare our method to two selected ones to show that our method is especially powerful to calculate the $\zeta$-values at large odd integer arguments.
\subsection*{Numerical test and error bound of the algorithm}
We regard the $\zeta$-values obtained by Mathematica as benchmarks. The result of the Riemann zeta function at odd integers with $n=1,2,\cdots,10$ obtained by our method(approximate value) is presented in table \ref{tab1}. The accuracy values and the absolute errors are also presented at the same time.
\begin{table}[h!]
\caption{Comparison between accurate values $\zeta^{\bf{ac}}(2n+1)$ and approximate values $\zeta^{\bf{ac}}(2n+1)$.}\label{tab1}
\begin{tabular}{|c|c|c|r|} 
\hline
$n$ &$\zeta^{\bf{ap}}(2n+1)$ &$\zeta^{\bf{ac}}(2n+1)$ &\multicolumn{1}{|c|}{Errors} \\
\hline
1	&1.201335874256	&1.202056903160	&-0.007210289040\\
\hline
2	&1.036972837734	&1.036927755143	&0.000045082590\\
\hline
3	&1.008365209797	&1.008349277382	&0.000015932415\\
\hline
4	&1.002011075857	&1.002008392826	&0.000002683031\\
\hline
5	&1.000494555053	&1.000494188604	&0.000000364486\\
\hline
6	&1.000122758824	&1.000122713348	&0.000000045476\\
\hline
7	&1.000030593607	&1.000030588236	&0.000000005371\\
\hline
8	&1.000007637815	&1.000007637198	&0.000000000617\\
\hline
9	&1.000001908283	&1.000001908213	&0.000000000070\\
\hline
10	&1.000000476941	&1.000000476933	&0.000000000008\\
\hline
\end{tabular}
\end{table}

Table \ref{tab1} tells us that, the idea that making the geometric mean instead of any of the upper or lower bound(see Theorem 2) be the best estimate of the Riemann zeta function dramatically reduces errors and satisfactory accuracy such as twelve decimal places in the tenth odd-argument of the Riemann zeta function can be achieved. It's interesting to find that only the Ap\'{e}ry's constant $\zeta(3)$ sightly larger than the approximate value obtained by our method. It's also funny to see the errors present an upside-down stair configuration, which implies that the error declines about ten times as long as the argument $n$ increase 1.

In table \ref{tab2} we present the absolute errors $ \epsilon(n)$ versus $n$, for the purpose of exploring the error bound when the argument $n$ is large enough.
\begin{table}[h!]
\caption{The errors of $\zeta(2n+1)$ based on our method.}\label{tab2}
\begin{tabular}{|l|l|l|l|}
\hline
\multicolumn{1}{|c|}{$n$}    &\multicolumn{1}{|c|}{Errors}   &\multicolumn{1}{|c|}{$n$}   &\multicolumn{1}{|c|}{Errors}\\
\hline
$1\times10^2$	&$1.05\times10^{-97}$	  &$1\times10^4$	&$1.04\times10^{-9544}$\\
\hline
$2\times10^2$	&$3.94\times10^{-193}$    &$2\times10^4$	&$3.92\times10^{-19087}$\\
\hline
$5\times10^2$	&$2.10\times10^{-479}$    &$5\times10^4$	&$2.08\times10^{-47714}$\\
\hline
$1\times10^3$	&$1.59\times10^{-956}$    &$1\times10^5$	&$1.56\times10^{-95426}$\\
\hline
$2\times10^3$	&$9.09\times10^{-1911}$   &$2\times10^5$	&$8.75\times10^{-190851}$\\
\hline
$5\times10^3$	&$1.70\times10^{-4773}$   &$5\times10^5$	&$1.55\times10^{-477123}$\\
\hline
\end{tabular}
\end{table}

It's clear that, from table \ref{tab2}, the error is of the order $O(10^{-n})$ approximately. By using of least square method, we notice that
\begin{equation}\label{eq:4errorbound}
\lg(\epsilon(n))=-0.9542n-1.6884.
\end{equation}
This formula suggests that when the argument of Riemann zeta function is large enough, our algorithm should be powerful enough to obtain the $\zeta$-values at odd integers.
\subsection*{Compare with the existed methods.}
In this subsection, we aim to compare our algorithm with the already existed ones, namely the Gauss-Hermite quadrature(Integral method, see \cite{integral}, Corollary 3.1) and rapid converging series(Series method, see \cite{series2}, eq.(3.30)). The Gauss-Hermite quadrature formula has the form\cite{integral}
\begin{equation}\label{GaussHermiteformula}
\int_{-\infty}^{\infty}f(x)e^{-x^2}{\rm{d}}x=\sum_{k=1}^{N}w_{k}f(x_k)+R_N
\end{equation}
where $x_k$ is one of the zeros of $H_N(x)$, the Hermite polynomial of degree $N$, and $w_k=-\frac{2^{N+1}N!\sqrt{\pi}}{H_N'(x_k)H_{N+1}(x_k)}$ is the corresponding weight. $R_N=\frac{N!\sqrt{\pi}}{2^{N}(2N)!}f^{2N}(\eta)$,$\eta\in(-\infty,\infty)$ is, obviously, the error bound of the above integral. Riemann zeta function is such an amazing function that it can be transformed into\cite{integral}
\begin{equation}\label{zetafuncInt}
\zeta(s)=\frac{\int_{-\infty}^{\infty}\big(\vert x\vert^{2s-1}{\rm{e}}^{-x^2}/(1-{\rm{e}}^{-x^2})\big){\rm{d}}x}{\int_{-\infty}^{\infty}\vert x\vert^{2s-1}{\rm{e}}^{-x^2}{\rm{d}}x}
\end{equation}
whose numerator and denominator are of the form presented in \eqref{GaussHermiteformula}. Among all the series representations of Riemann zeta function, the series below

\begin{align}\label{eqZetaSeries}
&\zeta(2n+1)=\frac{(-1)^{n-1}(2\pi)^{2n}}{(2n)![2^{2n}(2n-3)-2n+1]}\cdot\nonumber\\
&\Big[\sum_{m=1}^{n-1}(-1)^m\binom{2n-1}{2m-2}\frac{(2m)!(2^{2m}-1)}{(2\pi)^{2m}}\zeta(2m+1)+2\sum_{k=0}^{\infty}\frac{\zeta(2k)}{(2k+2n-1)(k+n)(2k+2n+1)2^{2k}}\Big]
\end{align}
converges most rapidly as pointed out by H.M. Srivastava\cite{series2}. When $n=1$ for instance, the error bound $R_N^{\bf(s)}$ of the $N$-th partial sum of the infinite series in \eqref{eqZetaSeries} satisfies
\begin{align}\label{errorboundSeries}
&\vert R_N^{\bf(s)}\vert =\frac{4\pi^2}{15}\sum_{k=N+1}^{\infty}\frac{\zeta(2k)}{(2k+1)(k+1)(2k+3)4^{k}}\nonumber\\
&                        <\frac{4\pi^2}{15}\frac{\zeta(2N+2)}{(2N+3)(N+2)(2N+5)}\sum_{k=N+1}^{\infty}\frac{1}{4^{k}}\nonumber\\
&                        =\frac{4\pi^2}{45}\frac{1}{(2N+3)(N+2)(2N+5)(4^{N}-\frac{1}{2})}
\end{align}
where we have used the fact that $\zeta(s)<\frac{1}{1-2^{1-s}}$ since $\eta(s)<1$. if $N=25$, the error bound is $\vert R_{25}^{\bf(s)}\vert < 1.0\times10^{-20}$, which is superior to other rapid series $\vert R_{25}\vert < 0.9\times10^{-18}$ as noted in \cite{series1,zeta3error}. Specially, when $N$ is larger than some typical numbers, the asymptotic behavior of \eqref{errorboundSeries} reads
\begin{equation}\label{errorboundSeriesAsym}
\lg(\vert R_N^{\bf(s)}\vert) \sim -2\lg2(N+1)-3\lg N.
\end{equation}
The accuracy of the latter two methods rely on the number of zeros(denoted as $N_1$) of the associated polynomial(in this occasion it is Hermite polynomial) and the terms(denoted as $N_2$) of partial sum of the infinite series respectively. We set the two integers be the same value, i.e. $N_1=N_2=25$ since the corresponding methods are both efficient as have been declared by many Mathematicians.
\begin{table}[h!]
\caption{Errors of three different methods for $\zeta(2n+1)$.}\label{tab3}
\begin{tabular}{|c|l|c|l|}
\hline
\multirow{2}{*}{$n$}
&\multirow{2}{*}{Integral method}
&Series method
&\multirow{2}{*}{Our method}\\
& &($\times10^{-20}$)&\\
\hline
3	&$2.42\times10^{-7}$	&$3.17434484$  &$1.50\times10^{-5}$\\
\hline
6	&$1.21\times10^{-10}$	&$3.14630746$  &$4.52\times10^{-8}$\\
\hline
9	&$4.31\times10^{-11}$	&$3.14592124$  &$4.99\times10^{-11}$\\
\hline
12	&$1.36\times10^{-12}$	&$3.14591532$  &$9.79\times10^{-14}$\\
\hline
15	&$7.40\times10^{-13}$	&$3.14591522$  &$1.35\times10^{-16}$\\
\hline
18	&$6.20\times10^{-13}$	&$3.14591522$  &$1.85\times10^{-19}$\\
\hline
21	&$8.46\times10^{-14}$	&$3.14591522$  &$2.54\times10^{-22}$\\
\hline
24	&$7.99\times10^{-15}$	&$3.14591522$  &$3.48\times10^{-25}$\\
\hline
27	&$7.77\times10^{-16}$	&$3.14591522$  &$4.78\times10^{-28}$\\
\hline
30	&$1.11\times10^{-16}$	&$3.14591522$  &$6.55\times10^{-31}$\\
\hline
\end{tabular}
\end{table}

The behaviors of the error bound of integral method and series method, as can be seen from table \ref{tab3}, are totally different. When the argument increases, the errors of the former decrease exponential from a high level, while the latter maintain at a nearly constant low level despite of the variation of $n$. Our method exhibits the worst results for small arguments, but the errors decrease dramatically with argument increasing. It outstrips integration method and series method before $n=12$ and $n=21$ respectively. Our method superior to them absolutely afterwards. To reach the accuracy obtained by our method, the number of nodes and terms in the above two methods should be augmented largely. In the series method for instance, the terms of the order $n$ in the infinity series should be included according to \eqref{eq:4errorbound} and \eqref{errorboundSeriesAsym}. Obviously, it is almost impossible to carry on within the limited CPU time when $n$ is an astronomical number.
\section*{Conclusion.}
In summary we firstly introduce two kinds of reduced Bernoulli numbers(RBNs) and prove their asymptotic behaviors in an uniform framework, and their series and integral representations are available at the same time. What's more, we discover and prove a recurrence formula \eqref{eq:4recur} of the Riemann zeta function original and construct an algorithm to evaluate the Riemann zeta function at odd integers based on it. The idea of our method is quiet simple, but it turns out to be a competent algorithm. The behavior of the error bound $\epsilon(n)$ is governed by $\lg(\epsilon(n))=-0.9542n-1.6884$ or $\epsilon(n)=O(10^{-n})$ approximately, which, of course, suggests that our method is especially suit for the calculation of $\zeta$-values at large odd integer arguments. Therefore, our results can also work as benchmarks to test the accuracy of other related algorithms. However, more works should be carried on to improve the accuracy at small arguments in future. Remarkably, the recurrence formula \eqref{eq:4recur} is likely to act as a touchstone to explore the closed form of the Riemann zeta function at positive integers since it witnesses the connection between $\zeta$-values at odd integers and even integers.

\section*{Acknowledgements}
The authors would like to show their appreciation to Junesang Choi, Yong Lin and Changle Liu for some useful discussions, and express their thanks to Jinlin Liu and Jiurong Han for their suggestions. Especially, they wishes to thank the anonymous referees of this paper for valuable suggestions which have improved the presentation of the paper.


\end{document}